\documentclass[16pt]{amsart}
\usepackage{graphicx,color}
\usepackage[left=1in,right=1in,top=1in,bottom=1in]{geometry}
\usepackage{latexsym}
\usepackage{cite}
\usepackage{amsmath,amssymb,amsthm,mathrsfs}

\newtheorem{example}{Example}

\allowdisplaybreaks[1]  

\newtheorem{theorem}{Theorem}

\newtheorem{remark}[theorem]{Remark}

\renewcommand{\phi}{\varphi}

\let\epsilon=\varepsilon

\newcommand{\Rc}{\operatorname{Rc}}
\newcommand{\Rm}{\operatorname{Rm}}

\def\crn#1#2{{\vcenter{\vbox{
\hbox{\kern#2pt \vrule width.#2pt height#1pt
 }
\hrule height.#2pt}}}}

\newcounter{mnotecount}[section]
\let\oldmarginpar\marginpar
\renewcommand\marginpar[1]{\-\oldmarginpar[\raggedleft\footnotesize #1]%
{\raggedright\footnotesize #1}}

\begin{document}

\title[Symbol of the RG-2 flow]{Short-time existence for the second order renormalization group flow in general dimensions.}
\author{Karsten Gimre}
\address[Gimre]{Department of Mathematics, Columbia University, New York City, New York}
\email{gimre@math.columbia.edu}
\author{Christine Guenther}
\address[Guenther]{Department of Mathematics and Computer Science, Pacific University,
Forest Grove, Oregon, 97116}
\email{guenther@pacificu.edu}
\author{James Isenberg}
\address[Isenberg]{Department of Mathematics, University of Oregon, Eugene, Oregon}
\email{isenberg@uoregon.edu}

\thanks{KG is partially supported by the NSF under grant DGE-1144155.}
\thanks{CG is partially supported by the Simons Foundation Collaboration Grant for Mathematicians - 283083}
\thanks{JI is partially supported by the NSF under grant PHY-1306441 at the University of Oregon. He also wishes to thank the Mathematical Sciences Research Institute in Berkeley, California, for support under grant 0932078 000. Some of this work was carried out while JI was in residence at MSRI during the fall of 2013.}

\date{\today}

\begin{abstract}
We prove local existence for the second order Renormalization Group flow initial value problem on closed Riemannian manifolds $(M,g)$ in general dimensions, for initial metrics whose sectional curvatures $K_P$ satisfy the condition $1+\alpha K_P > 0$, at all  points $p \in M$ and planes $P \subset T_p M$. This extends results previously proven for two and three dimensions. 
\end{abstract}

\maketitle

The second order approximation of the Renormalization Group flow for the nonlinear sigma model of quantum field theory, which we label  the RG-2 flow,  is specified  by the PDE system
\begin{equation}
\label{RG2Flow}
\frac{\partial }{\partial t}g=-2\Rc-\frac{\alpha }{%
2}\Rm^2. 
\end{equation}
Here $g$ is a Riemannian metric, $\Rc$ is its Ricci curvature,  $\Rm^2_{ij}=g^{pk}g^{ql}g^{nm}R_{iklm}R_{jpqn}$, and $\alpha$ is a positive parameter. We note that for our purposes here, $\alpha$ can assume any real value. For $\alpha =0$, this system \eqref{RG2Flow} reduces to the Ricci flow. One can see that the sign of the right hand side, which is roughly  $1 + \alpha \times$Curvature, should have an impact on the behavior of the flow, and this has been confirmed in various settings: in particular, the size of the term influences the parabolicity of the flow. Oliynyk has shown in \cite{O} that on a two-dimensional manifold,  if the Gaussian curvature $K$ satisfies the  condition $1+\alpha K > 0$, then the  flow is (weakly) parabolic; while  if  $1+\alpha K <0$ then  the flow is backward parabolic. In \cite{CM},  Cremaschi and Mantegazza prove that short-time existence holds in three dimensions so long as 
the analogous curvature condition $1+\alpha K_P > 0$ is satisfied for all sectional curvatures $K_P$. In this note we extend this curvature criterion for short-time existence for RG-2 flow  to all dimensions, as first announced in \cite{GGIa}. Our  main result is the following:  

\begin{theorem}
\label{ndim} 
Let $(M, g_0)$ be a closed $n$-dimensional Riemannian manifold. If  $1+\alpha K_P > 0$ for all sectional curvatures $K_P(g_0)$, at all points $p \in M$ and planes $P\subset T_p M$,  then there exists a unique solution $g(t)$ of the initial value problem $\frac{\partial }{\partial t}g=-2\Rc-\frac{\alpha }{2}\Rm^{2}$,  $g(0) = g_0$,  on some time interval $[0,T).$ 
\end{theorem}

\begin{remark}

In \cite{O}, Oliynyk finds open subspaces of the space of smooth metrics that are invariant under the two-dimensional RG-2 flow, and for which the flow remains parabolic (resp. backward parabolic). We are currently investigating this for general dimensions. 
\end{remark}

\begin{proof}

To prove the theorem, we calculate the principal symbol of the DeTurck-modified version of RG-2 flow, which is generated by the PDE system (compare with \eqref{RG2Flow} above)
\begin{equation}
\label{DeTurckRGF}
\frac{\partial }{\partial t}g_{ij}=-2R_{ij}+L_{W_{u,g}}g_{ij}-\frac{\alpha }{
2}\Rm_{ij}^2.  
\end{equation}%
Here $W_{u,g} = -g^{ij} u^{-1}_{jk}g^{kl}g^{pq}(\nabla_p u_{ql} - \frac{1}{2 }\nabla_l u_{pq})$ is the standard  vector field usually chosen to modify the Ricci flow into the related (parabolic) DeTurck version of Ricci flow, with $u$ a fixed metric.
Letting $\phi_t$ be the one-parameter family of diffeomorphisms 
generated by the vector field $-W_{u,g}$, then $\phi_t^* g$ is a solution of the RG-2 flow (see also \cite{GO}). As in the analogous Ricci flow case, if one can show (for a class of choices of the initial metric) that the PDE system \eqref{DeTurckRGF} is parabolic, then short-time existence holds for the RG-2 flow \eqref{RG2Flow}  as well as for the DeTurck-modified flow  \eqref{DeTurckRGF}.

To calculate the symbol of  the system \eqref{DeTurckRGF}, we first linearize the flow. For the first two terms of the right hand side of \eqref{DeTurckRGF}, this linearization effectively produces the Laplacian (see \cite{D}, or Theorem 2.1 in \cite{GIK}). For the remaining term, $\Rm^2$,  it is useful to recall the formula for the variation of  the Riemann curvature tensor with respect to the metric 
 (see pg.\ 74 in \cite{CK}): 
\begin{align}
\label{RmVar}
\nonumber  [D \Rm_g(h)]^l_{ijk} &=\left[\left.\frac{\partial}{\partial \epsilon}Rm(g+\epsilon h)\right|_{\epsilon = 0}\right]_{ijk}^l \\
 &=\frac{1}{2}g^{lp}(
\nabla _{i}\nabla _{j}h_{kp}+\nabla _{i}\nabla _{k}h_{jp}-\nabla _{i}\nabla
_{p}h_{jk}\\
&\notag\qquad-\nabla _{j}\nabla _{i}h_{kp}-\nabla _{j}\nabla _{k}h_{ip}+\nabla _{j}\nabla
_{p}h_{ik})
+ \text{LOT}.
\end{align} 
Here, we use $g$ to denote the metric about which we are linearizing and we use $h$ to denote the tangent to the linearization; we note that covariant derivatives and curvature terms appearing here and below are calculated with respect to $g$, and indices are raised and lowered using $g$. The term LOT denotes lower order terms with respect to derivatives of $h$. Applying \eqref{RmVar} together with the observations that  $[D g^{-1}( h)]^{ij} = -h^{ij}$ and   $\nabla_p\nabla_jh_n^l = \nabla_j\nabla_ph_n^l - R_{jpn}^m h_{m}^l  + R_{jpm}^l h_n^m$, we calculate 
\bigskip 
\begin{align}
\label{DRmSqr}
\nonumber [D \Rm^2_g(h)] _{ij} &=\frac{1}{2}g^{pk}R_{ikl}^{n}(\nabla _{j}\nabla _{p}h_{n}^{l}+\nabla
_{j}\nabla _{n}h_{p}^{l}-\nabla _{j}\nabla ^{l}h_{pn}\\
\nonumber &\qquad -\nabla _{p}\nabla _{j}h_{n}^{l}-\nabla _{p}\nabla _{n}h_{j}^{l}+\nabla
_{p}\nabla ^{l}h_{jn}) \\
\nonumber&\qquad+\frac{1}{2}
g^{pk}R_{jpn}^{l}(\nabla _{i}\nabla _{k}h_{l}^{n}+\nabla _{i}\nabla
_{l}h_{k}^{n}-\nabla _{i}\nabla ^{n}h_{kl} \\
\nonumber&\qquad-\nabla _{k}\nabla _{i}h_{l}^{n}-\nabla _{k}\nabla _{l}h_{i}^{n}+\nabla
_{k}\nabla ^{n}h_{il}) + \text{LOT} \\
&= R_{iklu}(\nabla_j\nabla^l h^{ku} - \nabla^k \nabla^l h_j^u) + R_{jklu}(\nabla_i\nabla^l h^{ku} - \nabla^k\nabla^lh_i^u) +\text{LOT}.
\end{align}%

 We obtain the principal symbol of the full flow \eqref{DeTurckRGF} by replacing  each $\nabla$ appearing in the sum of the Laplacian plus the  terms in (\ref{DRmSqr}) by  the co-vector $\xi$. As in the Ricci flow case (see \cite{H}), we work in orthonormal coordinates,  and without loss of generality, we assume that $\xi_1 = 1, \ \xi_l = 0, \ l \ge 2$. Writing the derivative of the right hand side of (\ref{DeTurckRGF}) as $DL_g$, the symbol is thus given by
 \begin{equation}
 \label{sdf}
\sigma DL_g(h)_{ij} :=\sigma DL_g(\xi)(h)_{ij}  = h_{ij}+\frac{\alpha}{2} R_{ik1u}\delta_{j1}h^{ku}-\frac{\alpha}{2} R_{i11u}h_j^u+  \frac{\alpha}{2} R_{jk1u}\delta_{i1}h^{ku} - \frac{\alpha}{2} R_{j11u}h_i^u.
\end{equation}
Noting that (by design) the left hand side of (\ref{sdf}) takes the form of a linear algebraic operator on $h$, we rewrite the right hand  side of \eqref{sdf} formally as a matrix expression, 
\begin{equation} 
\label{matrix}
\sigma DL_g(h)_A=\Sigma_A^B h_A.
\end{equation}
The capital Latin letters represent the symmetric indices on $h$:
\begin{equation}
h_A \leftrightarrow \{ h_{11}, h_{12}, ...,h_{1n}, h_{22}, h_{23}, ..., h_{2n},...,h_{nn}\},
\end{equation}
with, for example, $h_{21}$ not appearing since $h_{12}$ does.  We note that $h_A $ is an $n(n+1)/2$ dimensional vector, and correspondingly $\Sigma_A^B$ is an $n(n+1)/2 \times n(n+1)/2$ matrix, whose columns are given by the $h_{B}$ terms of $\sigma DL_g(h)_A$. 

To determine whether or not the PDE system \eqref{DeTurckRGF} is parabolic for a given metric $g$ (representing the initial data for an RG-2 flow solution) one determines if the matrix $\Sigma_A^B $ is nondegenerate for that metric. We verify  here that the condition stated for sectional curvatures in the hypothesis of Theorem \ref{ndim} guarantees this nondegeneracy, and consequently the parabolicity of \eqref{DeTurckRGF}. To carry out this verification explicitly, it is useful to write out the components of equation (\ref{sdf}) for various possible choices of the indices $ij$. There are three cases to consider:

\bigskip

Case 1: $i = j = 1$:

\begin{equation}
\sigma DL_g(h)_{11} = h_{11}+ \alpha R_{1k1u}h^{ku},
\end{equation}

Case 2: $ i = 1$, $j \ne 1$:

\begin{equation}
\sigma DL_g(h)_{1j}=h_{1j} + \frac{\alpha}{2}R_{jk1u}h^{ku},
\end{equation}
where $k \ge 2,$ 

Case 3: $ i \ne 1$, $j \ne 1$

\begin{equation}
\label{Case3}
\sigma DL_g(h)_{ij} = h_{ij} + \frac{\alpha}{2} R_{1i1u} h_j^u+ \frac{\alpha}{2} R_{1j1u} h_i^u .
\end{equation}

Based on these expressions, one sees that the matrix $\Sigma_A^B$ takes the block form
\begin{equation}
\label{matrixout}
\Sigma_A^B=\begin{pmatrix}I&\lambda\\ \mu&\nu\end{pmatrix}.
\end{equation}
Here  $I$ is the $n \times n$ identity matrix, and since  $n(n+1)/2 - n = n(n-1)/2$, the block $\mu$ is an $n \times n(n-1)/2$ matrix composed entirely  of  zeros. Thus the determinant of $\Sigma_A^B$
is given by $\det(\nu)$, from which it follows that the system \eqref{DeTurckRGF} is parabolic so long as the $ n(n-1)/2 \times n(n-1)/2$ matrix $\nu$  has nonvanishing determinant. 
One determines the components  of $\nu$  from expression \eqref{Case3}. The diagonal components of $\nu$ have $A=B={ij}$ (so they are the $h_{ij}$ terms of $\sigma DL_g(h)_{ij}$ for $i,j >n$), and have the form
\begin{equation}
 1 + \frac{\alpha}{2} R_{1i1i}+ \frac{\alpha}{2} R_{1j1j}. 
\end{equation}
These terms thus take the form $1 + \alpha K_A$, for sectional curvatures $K_A$.

By contrast, the off-diagonal terms of $\nu$ do not involve the sectional curvatures; rather, they involve only mixed curvatures of the form $R_{1i1u}$ and $R_{1j1u}$ for $u\neq i$ and $u \ne j$.
Standard arguments show that we can choose an orthonormal basis which diagonalizes the (symmetric) matrix  $R_{1m1n}$; it follows that $\nu$ is nondegenerate so long as the  condition $1 + \alpha K_A>0$ holds, thus proving our theorem.



\end{proof}

\begin{example}{4-Dimensions}

To help clarify the discussion in the proof above, we write out the matrix $\nu$ explicitly for  n = 4 dimensions. The order of the columns here is $h_{22}, h_{23}, h_{24}, h_{33}, h_{34}, h_{44}$. Thus, for  example,  the first row in this matrix is obtained by writing out $ \sigma DL_g(h)_{22} = (1+\alpha  R_{1212}) h_{22} +  \alpha R_{1213}h_{23} +  \alpha R_{1214}h_{24}.$ The complete set of entries are as follows: 

$$\nu=\begin{pmatrix}
1+\alpha R_{1212} & \alpha R_{1213} & \alpha R_{1214} &0 &0 &0\\
\frac{\alpha}{2} R_{1312} & 1+\frac{\alpha}{2} (R_{1313} + R_{1212}) & \frac{\alpha}{2}  R_{1314} & \frac{\alpha}{2} R_{1213} & \frac{\alpha}{2} R_{1214} & 0 \\
\frac{\alpha}{2} R_{1214} & \frac{\alpha}{2} R_{1314} & 1+\frac{\alpha}{2} (R_{1212} + R_{1414}) & 0 & \frac{\alpha}{2} R_{1213} & \frac{\alpha}{2} R_{1214} \\
0&\alpha R_{1213} & 0 & 1+\alpha R_{1313} & \alpha R_{1314} & 0 \\
0 & \frac{\alpha}{2} R_{1214} & \frac{\alpha}{2} R_{1213} & \frac{\alpha}{2} R_{1314} & 1+\frac{\alpha}{2} (R_{1313} + R_{1414}) & \frac{\alpha}{2} R_{1314} \\
0 & 0 & \alpha R_{1214} & 0 &  \alpha R_{1314} & 1+ \alpha R_{1414} 
\end{pmatrix}$$
We can diagonalize $R_{1m1n}$, ensuring that the off-diagonal terms are zero. 

\end{example}

\end{document}